\documentclass{article}
\usepackage{amsmath,amssymb,graphics,graphicx}
\newtheorem{df}{Definition}[section]
\newtheorem{thm}[df]{Theorem}
\newtheorem{lem}[df]{Lemma}
\newtheorem{prop}[df]{Proposition}
\newtheorem{cor}[df]{Corollary}
\newcommand{\cAa}{\textbf{\textsf{A1}}}
\newcommand{\cAb}{\textbf{\textsf{A2}}}
\newcommand{\cAc}{\textbf{\textsf{A3}}}
\newcommand{\cAd}{\textbf{\textsf{A4}}}
\newcommand{\cBa}{\textbf{\textsf{B1}}}
\newcommand{\cBb}{\textbf{\textsf{B2}}}
\newcommand{\cBc}{\textbf{\textsf{B3}}}
\newcommand{\cBd}{\textbf{\textsf{B4}}}
\newcommand{\cCa}{\textbf{\textsf{C1}}}

\newcommand{\cCd}{\textbf{\textsf{C4}}}
\newcommand{\cDa}{\textbf{\textsf{D1}}}

\newcommand{\cDd}{\textbf{\textsf{D4}}}

\title{Semi-continuous and discrete wavelet frames on $n$-dimensional spheres}
\author{I. Iglewska-Nowak\footnote{West Pomeranian University of Technology in Szczecin, School of Mathematics, al. Piast\'ow 17, 70--310 Szczecin, Poland}}

\begin{document}

\maketitle

\bibliographystyle{amsplain}

\begin{abstract} The paper shows that under some mild conditions $n$-dimensional spherical wavelets derived from approximate identities build semi-continuous frames. Moreover, for sufficiently dense grids Poisson wavelets on $n$-dimensional spheres constitute a discrete frame. In the proof we only use the localization properties of the reproducing kernel and its gradient.
\end{abstract}

\begin{bfseries}Keywords:\end{bfseries} spherical wavelets, Poisson wavelets, frames, $n$-spheres \\
\begin{bfseries}AMS Classification:\end{bfseries} 42C40

\section{Introduction}
In the recent years, an interest on $n$--dimensional spherical wavelet transform has been growing. Besides discrete approaches \cite{NW96,IINMul} there are several continuous constructions \cite{AVn} (being a generalization to $n$ dimensions of spherical wavelets introduced in~\cite{AV}), \cite{sE11,EBCK09}. For an efficient usage of a continuous wavelet transform, a discretization algorithm is needed. Frames have been constructed for $2$--dimensional spherical wavelets derived in~\cite{AV}, cf. \cite{jpA,BVAJM}, however, the phase--space discretization is performed on an equiangular grid, a solution that can hardly be applied in a higher dimension.

In this paper, we generalize the results obtained in \cite{IH10} for $2$--dimensional spherical wavelets. The construction of semi--continuous frames is similar to that in~\cite{jpA,BVAJM} for the two--dimensional sphere. As a next step, for each scale we perform a discretization of the spherical parameter such that the sampling points are quite uniformly distributed over the sphere. Finally, the sampling point positions are perturbed in such a way that the density of the resulting grid is controlled with respect to the scale and space parameter simultaneously. If the density is big enough, the discrete set of wavelets is a frame for~$\mathcal L^2(\mathcal S^n)$. The constraints on the wavelets are some estimations on their reproducing kernel and its gradient, which are satisfied by Poisson multipole wavelets.

After the present research had been completed, the author learnt about a similar frame construction for Mexican needlets~\cite{GM09b}. Therefore, the present discussion focuses on the differences and similarities in the end of the paper.

The paper is organized as follows: In Section~\ref{sec:sphere} after a recapitulation of basic facts about the $n$--dimensional wavelet transform derived from approximate identities, and particularly Poisson multipole wavelets, we recall some information about frames in Hilbert spaces. Section~\ref{sec:semicontinuous_frames} contains a discussion of a condition for semi--continuous frames. It is shown, that many popular wavelet families constitute semi--continuous frames. The main theorem of this paper about the phase--space discretization is to be found in Section~\ref{sec:discrete_frames}, and a perturbation of this result to fully irregular frames controlled only by hyperbolic density is the topic of Section~\ref{sec:density}. It is shown in Section~\ref{sec:Poisson_frames} that both discretization results apply to Poisson multipole wavelets. A comparison with the frame construction for Mexican needlets is presented in Section~\ref{sec:Mexican_needlets}.

\section{Preliminaries}\label{sec:sphere}

By $\mathcal{S}^n$ we denote the $n$--dimensional unit sphere in $(n+1)$--dimensional Euclidean space~$\mathbb{R}^{n+1}$ with the rotation--invariant measure~$d\sigma$ normalized such that
$$
\Sigma_n=\int_{\mathcal{S}^n}d\sigma=\frac{2\pi^{\lambda+1}}{\Gamma(\lambda+1)},
$$
where~$\lambda$ and~$n$ are related by
$$
\lambda=\frac{n-1}{2}.
$$
The surface element $d\sigma$ is explicitly given by
$$
d\sigma=\sin^{n-1}\theta_1\,\sin^{n-2}\theta_2\dots\sin\theta_{n-1}d\theta_1\,d\theta_2\dots d\theta_{n-1}d\varphi,
$$
where $(\theta_1,\theta_2,\dots,\theta_{n-1},\varphi)\in[0,\pi]^{n-1}\times[0,2\pi)$ are spherical coordinates satisfying
\begin{align*}
x_1&=\cos\theta_1,\\
x_2&=\sin\theta_1\cos\theta_2,\\
x_3&=\sin\theta_1\sin\theta_2\cos\theta_3,\\
&\dots\\
x_{n-1}&=\sin\theta_1\sin\theta_2\dots\sin\theta_{n-2}\cos\theta_{n-1},\\
x_n&=\sin\theta_1\sin\theta_2\dots\sin\theta_{n-2}\sin\theta_{n-1}\cos\varphi,\\
x_{n+1}&=\sin\theta_1\sin\theta_2\dots\sin\theta_{n-2}\sin\theta_{n-1}\sin\varphi.
\end{align*}
$\left<x,y\right>$ or $x\cdot y$ stand for the scalar product of vectors with origin in~$O$ and endpoints on the sphere. As long as it does not lead to misunderstandings, we identify these vectors with points on the sphere. By $\angle(x,y)$ we denote the geodesic distance between two points $x,y\in\mathcal S^n$,
$$
\angle(x,y):=\arccos(x,y).
$$

Scalar product of $f,g\in\mathcal L^2(\mathcal S^n)$ is defined by
$$
\left<f,g\right>_{\mathcal L^2(\mathcal S^n)}=\frac{1}{\Sigma_n}\int_{\mathcal S^n}\overline{f(x)}\,g(x)\,d\sigma(x),
$$
and by $\|\circ\|$ we denote the induced $\mathcal L^2$--norm.

Gegenbauer polynomials $C_l^\lambda$ of order~$\lambda\in\mathbb R$ and degree~$l\in\mathbb{N}_0$, are defined in terms of their generating function
$$
\sum_{l=0}^\infty C_l^\lambda(t)\,r^l=\frac{1}{(1-2tr+r^2)^\lambda},\qquad t\in[-1,1].
$$
They are real-valued and for some fixed $\lambda\ne0$ orthogonal to each other with respect to the weight function~$\left(1-\circ^2\right)^{\lambda-\frac{1}{2}}$, compare~\cite[formula~8.939.8]{GR}.

Let $Q_l$ denote a polynomial on~$\mathbb{R}^{n+1}$ homogeneous of degree~$l$, i.e., such that $Q_l(az)=a^lQ_l(z)$ for all $a\in\mathbb R$ and $z\in\mathbb R^{n+1}$, and harmonic in~$\mathbb{R}^{n+1}$, i.e., satisfying $\Delta Q_l(z)=0$, then $Y_l(x)=Q_l(x)$, $x\in\mathcal S^n$, is called a hyperspherical harmonic of degree~$l$. The set of hyperspherical harmonics of degree~$l$ restricted to~$\mathcal S^n$ is denoted by $\mathcal H_l=\mathcal H_l(\mathcal S^n)$. $\mathcal H_l$--functions are eigenfunctions of the Laplace--Beltrami operator $\Delta^\ast:=\left.\Delta\right|_{\mathcal S^n}$ with eigenvalue $-l(l+2\lambda)$, further, hyperspherical harmonics of distinct degrees are orthogonal to each other. The number of linearly independent hyperspherical harmonics of degree~$l$ is equal to
$$
N=N(n,l)=\frac{(n+2l-1)(n+l-2)!}{(n-1)!\,l!}.
$$

The addition theorem states that
$$
C_l^\lambda(x\cdot y)=\frac{\lambda}{l+\lambda}\,\sum_{\kappa=1}^N\overline{Y_l^\kappa(x)}\,Y_l^\kappa(y),
$$
for any orthonormal set $\{Y_l^\kappa\}_{\kappa=1,2,\dots,N(n,l)}$ of hyperspherical harmonics of degree~$l$ on~$\mathcal S^n$.
In this paper, we will be working with the orthonormal basis for~$\mathcal L^2(\mathcal S^n)=\overline{\bigcup_{l=0}^\infty\mathcal H_l}$, consisting of hyperspherical harmonics given by
$$
Y_l^k(x)=A_l^k\prod_{\nu=1}^{n-1}C_{k_{\nu-1}-k_\nu}^{\frac{n-\nu}{2}+k_\nu}(\cos\theta_\nu)\sin^{k_\nu}\theta_\nu\cdot e^{\pm ik_{n-1}\varphi}
$$
with $l=k_0\geq k_1\geq\dots\geq k_{n-1}\geq0$, $k$ being a sequence $(k_1,\dots,\pm k_{n-1})$ of integer numbers, and normalization constants~$A_l^k$, compare~\cite{Vilenkin,AI}.

Every $\mathcal{L}^1(\mathcal S^n)$--function~$f$ can be expanded into Laplace series of hyperspherical harmonics by
$$
f\sim\sum_{l=0}^\infty f_l,
$$
where $f_l$ is given by
$$
f_l(x)=\frac{\Gamma(\lambda)(l+\lambda)}{2\pi^{\lambda+1}}\int_{\mathcal S^n}C_l^\lambda(x\cdot y)\,f(y)\,d\sigma(y)
   =\frac{l+\lambda}{\lambda}\left<C_l^\lambda(x\cdot\circ),f\right>.
$$
For zonal functions (i.e., those depending only on~$\theta_1=\left<\hat e,x\right>$, where~$\hat e$ is the North Pole of the sphere
$\hat e=(1,0,\dots,0)$) we obtain the Gegenbauer expansion
$$
f(\cos\theta_1)=\sum_{l=0}^\infty\widehat f(l)\,C_l^\lambda(\cos\theta_1)
$$
with Gegenbauer coefficients
$$
\widehat f(l)=c(l,\lambda)\int_{-1}^1 f(t)\,C_l^\lambda(t)\left(1-t^2\right)^{\lambda-1/2}dt,
$$
where~$c$ is a constant that depends on~$l$ and~$\lambda$.

For $f,g\in\mathcal L^1(\mathcal S^n)$, $g$ zonal, their convolution $f\ast g$ is defined by
\begin{equation*}
(f\ast g)(x)=\frac{1}{\Sigma_n}\int_{\mathcal S^n}f(y)\,g(x\cdot y)\,d\sigma(y).
\end{equation*}
With this notation we have
$$
f_l(x)=\frac{l+\lambda}{\lambda}\,\left(f\ast C_l^\lambda\right)(x),
$$
i.e., the function $\mathcal K_l^\lambda:=\frac{l+\lambda}{\lambda}\,C_l^\lambda$ is the reproducing kernel for~$\mathcal H_l$.

Further, any function $f\in\mathcal L^2(\mathcal S^n)$ has a unique representation as a mean--convergent series
$$
f(x)=\sum_k a_l^kY_l^k(x),\qquad x\in\mathcal S^n,
$$
where
$$
a_l^k=a_l^k(f)=\frac{1}{\Sigma_n}\int_{\mathcal S^n}\overline{Y_l^k(x)}\,f(x)\,d\sigma(x)=\left<Y_l^k,f\right>,
$$
for proof cf.~\cite{Vilenkin}. In analogy to the two-dimensional case, we call $a_l^k$ the Fourier coefficients of the function~$f$. The Parseval identity has the form
\begin{equation}\label{eq:Parseval}
\|f\|^2_{\mathcal L^2(\mathcal S^n)}=\sum_{l=0}^\infty\sum_{\kappa=1}^{N(n,l)}|\alpha_l^\kappa|^2
\end{equation}

We identify zonal functions with functions over the interval $[-1,1]$, i.e., whenever it does not lead to mistakes, we write
$$
f(x)=f(\cos\theta_1).
$$

For further details on this topic we refer to the textbooks~\cite{Vilenkin} and~\cite{AH12}.

\subsection{Continuous wavelet transform on $n$--spheres}
Spherical wavelet transform derived from singular integrals was introduced by Freeden and Windheuser in~\cite{FW96} and~\cite{FW-C} (compare also~\cite{FGS-book}) for two--dimensional spheres and by Bernstein in~\cite{sB09} for three--dimensional spheres. A generalization to non--zonal wavelets in $n$ dimensions is to be found in~\cite{EBCK09} and in~\cite{IIN14CWT}, and we refer to these papers for more details.

We present here the basic facts for zonal wavelets with respect to the weight function $\alpha(a)=\frac{1}{a}$.

\begin{df}\label{def:sphwavelet}
A subfamily $\{\Psi_a\}_{a\in\mathbb{R}_+}$ of the space of zonal integrable functions is called a spherical wavelet if it satisfies the following admissibility conditions:
\begin{enumerate}
\item for $l\in\mathbb{N}_0$
\begin{equation}\label{eq:admwv1}
\int_0^\infty\bigl|\widehat{\Psi_a}(l)\bigr|^2\,\frac{da}{a}=\left(\frac{l+\lambda}{\lambda}\right)^2,
\end{equation}
\item for $R\in\mathbb{R}_+$
\begin{equation*}
\int_{-1}^1\left|\int_R^\infty\bigl(\overline{\Psi_a}\ast\Psi_a\bigr)(t)\,\frac{da}{a}\right|\,\left(1-t^2\right)^{\lambda-1/2}dt\leq M
\end{equation*}
with $M$ independent from~$R$. The spherical wavelet transform
\begin{equation}\label{eq:wvtrafo}
\mathcal{W}_\Psi:\,\mathcal{L}^2\left(\mathcal{S}^n\right)\to\mathcal{L}^2\left(\mathbb{R}_+\times\mathcal{S}^n\right)
\end{equation}
is defined by
$$
\mathcal{W}_\Psi f(a,x)=\left<\Psi_{a,x},f\right>_{\mathcal{L}^2(\mathcal{S}^n)}=\left(f\ast\overline{\Psi_a}\right)(x).
$$
\end{enumerate}
\end{df}

Examples of spherical wavelets are Abel--Poisson wavelets with Gegenbauer coefficients of the kernel
$$
\widehat{\Psi_a}(l)=\frac{l+\lambda}{\lambda}\,\sqrt{2la}\,e^{-la},\qquad l\in\mathbb N_0,
$$
Gauss--Weierstrass wavelets given by
$$
\widehat{\Psi_a}(l)=\frac{l+\lambda}{\lambda}\,\sqrt{2l(l+2\lambda)a}\,e^{-l(l+2\lambda)a},\qquad l\in\mathbb N_0,
$$
and Poisson multipole wavelets.

\subsection{Poisson wavelets on $n$--spheres}
Poisson multipole wavelets on two--dimensional spheres have proven to be very useful in applications~\cite{HCM03,CPMDHJ05}. Therefore, the author of the present paper defined Poisson wavelets on $n$--spheres in~\cite{IIN14PW} in a similar way as in~\cite{HI07}, i.e., as derivatives of Poisson kernel. Here, we recall the definition and some of their properties, which will be used for the proof of the main theorem. The statements come from~\cite{IIN14PW} and their proofs are similar to those for the two--dimensional case, given in~\cite{IH10}.

\begin{df}The Poisson wavelet of order~$m$, $m\in\mathbb N$, at a scale~$a$, $a\in\mathbb R_+$, is given recursively by
\begin{align*}
g_a^1&=a r\partial_ rp_{ r\hat e},\qquad r=e^{-a},\\
g_a^{m+1}&=a r\partial_ rg_a^m,
\end{align*}
where
$$
p_{ r\hat e}(y)=\frac{1}{\Sigma_n}\frac{1- r^2}{| r\hat e-y|^{n+1}}
   =\frac{1}{\Sigma_n}\frac{1- r^2}{(1-2 r\cos\theta+ r^2)^{(n+1)/2}},\quad\theta=\left<\hat e,y\right>.
$$
\end{df}

\begin{lem}The Gegenbauer expansion of Poisson wavelets is given by
\begin{align*}
g_a^{m}(y)&=\frac{1}{\Sigma_n}\,\sum_{l=0}^\infty\frac{l+\lambda}{\lambda}\,(al)^m e^{-al}\,C_l^\lambda(\cos\theta)\\
&=\frac{1}{\Sigma_n}\,\sum_{l=0}^\infty(al)^m e^{-al}\,\mathcal K_l^\lambda(\cos\theta).
\end{align*}
\end{lem}

Although a normalization constant is required for such wavelets to satisfy condition~\eqref{eq:admwv1}, we prefer to use formula~\eqref{eq:wvtrafo} with~$g_a^m$ defined as above and include the factor in the reproducing formula:
$$
f(x)=\frac{4^m\,\Sigma_n}{\Gamma(2m)}\int_0^\infty\!\!\int_{\mathcal S^n}\mathcal W_mf(a,y)\,g_{a,y}^m(x)\,\frac{d\sigma(y)\,da}{a}
$$
with $\mathcal W_m$ denoting $\mathcal W_{g^m}$. The image of the wavelet transform with respect to Poisson wavelets is  a Hilbert space with the reproducing kernel
\begin{equation}\label{eq:kernel}
\Pi^m(a,x;b,y)=\frac{4^m\,\Sigma_n^2}{\Gamma(2m)}\cdot\left<g_{a,x}^m,g_{b,y}^m\right>=\frac{4^m\,\Sigma_n}{\Gamma(2m)}\frac{(ab)^m}{(a+b)^{2m}}\,g_{a+b}^{2m}(x\cdot y).
\end{equation}
Further, the wavelets are space localized in the following way.

\begin{thm}\label{thm:scaling}
Let~$g_a^m$ be a Poisson wavelet of order~$m$. Then there exists a constant~$\mathfrak{c}$ such that
$$
|a^ng_a^m\left(\cos(a\theta)\right)|\leq\frac{\mathfrak{c}\cdot e^{-a}}{\theta^{m+n}},\quad\theta\in\left(0,\frac{\pi}{a}\right],
$$
uniformly in~$a$. $m+n$ is the largest possible exponent in this inequality.
\end{thm}

\begin{cor}\label{cor:proporcja} The functions $(a,\theta)\mapsto a^ng_a^m(\cos\theta)$ are bounded by $\mathfrak c\cdot e^{-a}$ uniformly in~$\theta$, and $n$ is the smallest possible exponent in this estimation.
\end{cor}

\subsection{Frames in reproducing kernel Hilbert spaces}
For our considerations we need the following characterization of general frames in reproducing kernel Hilbert spaces. The statements in this section come from~\cite{IH10}.

\begin{df} A family of vectors $\{g_x,\,x\in X\}\subset\mathcal{H}$ in a Hilbert space $\mathcal{H}$ indexed by a measure space $X$ with a positive measure $\mu$
is called a~\emph{frame with weight}~$\mu$ if the mapping $x\mapsto g_x$ is weakly measurable, i.e., $x\mapsto\left<g_x,u\right>$ is measurable, and if
for some $0\leq \epsilon<1$ we have
\begin{equation}  \label{eq:framedef}
(1-\epsilon)\,\|u\|^2\leq \int_{X}|\left<g_x,u\right>|^2 d\,\mu(x) \leq(1+\epsilon)\,\|u\|^2.
\end{equation}
for all $u\in\mathcal{H}$. Equivalently, the frame condition reads
\begin{equation*}
\left| \int_X |\left<g_x,u\right>|^2 \,d\mu(x) - \|u\|^2 \right| \leq\epsilon\, \|u\|^2.
\end{equation*}
If $\epsilon=0$, we call it a \emph{tight frame}. The numbers $A:=1-\epsilon$ and $B:=1+\epsilon$ are called \emph{frame bounds}.
\end{df}

\begin{bfseries}Remark.\end{bfseries} Usually in the frame theory frame bounds are only supposed to be positive numbers, cf. e.g.~\cite{CJ99}. In the case of weighted frames the weight can be scaled in such a way that the frame bounds satisfy additionally $A<1<B$.

Let $\mathcal H=\mathcal L^2(X,d\mu)$ be a Hilbert space of functions over $X$ with the reproducing kernel $\Pi$
\begin{equation*}
u(x) = \int_X \Pi(x,y)\,u(y)\,d\mu(y).
\end{equation*}
The family of functions $\{g_x = \Pi(x,\cdot)\}$ with $x\in X$ is a tight frame with weight~$\mu$. Conversely, a tight frame $\{g_x,\,x\in X\}$ and a measure $\mu$ in a Hilbert space~$\mathcal{H}$ are naturally associated with a reproducing kernel Hilbert space of functions in $\mathcal L^2(X,d\mu)$, as shown by the next theorem.

\begin{thm}
The mapping
\begin{equation}\label{eq:mappingF}
\mathcal{F} : \mathcal{H} \to \mathcal L^2(X,d\mu),\qquad\mathcal{F}u(x)=\left<g_x,u\right>
\end{equation}
is a partial isometry and the image $\mathcal{U}$ of this mapping is characterized by the reproducing kernel
\begin{equation*}
\Pi(x,y) = \left<g_x, g_y\right>.
\end{equation*}
That means, $u\in \mathcal L^2(X,d\mu)$ is in the range of $\mathcal{F}$ if and only if
\begin{equation*}
\int_X \Pi(x,y)\, u(y)\, d\mu(y) = u(x).
\end{equation*}
\end{thm}

The last integral is absolutely convergent since $\Pi(x,\cdot)$ is in $\mathcal L^2(X,d\mu)$.

In particular, we have:

\begin{prop}
Let $\{g_x,\,x\in X\}$, be a tight frame with weight $\mu$ on $\mathcal{H}$. A family $\{g_y\}$ with $y\in\Lambda\subset X$ and a measure $\lambda$ on $\Lambda$ yield
a frame for $\mathcal{H}$ if and only if $\{\Pi(y,\cdot),\,y\in \Lambda\}$, $\Pi(\xi,\eta)=\left<g_\xi,g_\eta\right>$, is a frame for $\mathcal{F}(\mathcal{H})$, with $\mathcal{F}$ given by~\eqref{eq:mappingF}.
\end{prop}

Frames of the form $\{\Pi(y,\cdot)\}$ can be characterized as follows:

\begin{thm}
Let $\Lambda\subset X$ and let $\lambda$ be a measure on $\Lambda$, and $\mu$ be a measure on~$X$. The family of functions $\{g_y=\Pi(y, \cdot),\,y\in\Lambda\}\subset \mathcal L^2(X,d\mu)$, is a frame with weight $\lambda$ for $\mathcal U=\mathcal F(\mathcal H)$ if and only if
\begin{equation}  \label{fff}
F(x,z) = \int_\Lambda \Pi(x,y)\,\Pi(y,z)\, d\lambda(y) - \Pi(x,z),
\end{equation}
is the kernel of a bounded operator~$\mathbb{F}$ on $\mathcal{U}$ with~$\|\mathbb{F}\|<1$.
\end{thm}

Since $\Pi(x,z) = \int_X \Pi(x,y)\,\Pi(y,z)\, d\mu(y)$, the theorem shows that the existence of frames is intimately linked to the
existence of good quadrature rules for functions in $\mathcal{U}$. This general principle will be used together with the following perturbation result.

\begin{cor}\label{cor:discretization} Suppose, for a set $\Lambda$ the family $\{g_y=\Pi(y, \cdot),\,y\in \Lambda\}$, is a weighted frame for $\mathcal{U}$ with weight $\lambda$. If now for another set $\Upsilon$ we have for $\{g_y=\Pi(y,\cdot),\,y\in\Upsilon\} \subset\mathcal{U}$ and a weight $\upsilon$ that
\begin{equation*}
G(x,z) = \int_{\Lambda} \Pi(x,y)\, \Pi(y,z)\, d\lambda(y) - \int_{\Upsilon} \Pi(x,y)\, \Pi(y,z)\,d\upsilon(y)
\end{equation*}
is the kernel of an operator~$\mathbb{G}$ with operator norm $\|\mathbb{G} \|\leq 1-\|\mathbb{F}\|$, where the kernel of $\mathbb{F}$ is given by~\eqref{fff}, then $\{g_y,\,y\in\Upsilon\}$ is a frame with weight $\upsilon$.
\end{cor}

More details on this topic can be found in~\cite{IH10}, \cite{CJ99}, and~\cite{oC03}.

\section{Semi--continuous frames}\label{sec:semicontinuous_frames}

In this section we prove the existence of semi--continuous wavelet frames, i.e., such that only the scale parameter is discretized.

\begin{lem}Let $\{\Psi_a\}$ be a wavelet family such that
\begin{equation}\label{eq:pseudo_generating}
\widehat\Psi_a(l)=\frac{l+\lambda}{\lambda}\cdot\gamma\bigl(a\cdot\tau(l)\bigr)
\end{equation}
for an arbitrary function~$\tau$, and let $\mathcal B=\{b_j\}_{j\in\mathbb N_0}$ be a countable set of scales. The set
\begin{equation}\label{eq:set_semifr}
J=\{\Psi_{b_j}(x)\colon b_j\in\mathcal B,\,x\in\mathcal S^n\}
\end{equation}
is a semi--continuous frame with weights~$\nu_j$ and bounds~$A$, $B$ if and only if
\begin{equation}\label{eq:semiframe_cond}
A\leq\sum_{j=0}^\infty\left|\gamma\bigl(b_j\tau(l)\bigr)\right|^2\,\nu_j\leq B
\end{equation}
holds independently of~$l$.
\end{lem}

This lemma (and its proof which we write here for the convenience of the reader) is a slight modification of \cite[Proposition~6]{IIN07} that has also been implicitly used in the proof of~\cite[Theorem 4]{IH10}.

\begin{bfseries}Proof.\end{bfseries}  Suppose, $J$ is a semi--continuous frame with bounds~$A$ and~$B$. Then the following holds
\begin{equation}\label{eq:semifr_cond}
A\|f\|^2\leq\sum_{j=0}^\infty\int_{\mathcal S^n}|\mathcal W_\Psi f(b_j,x)|^2\,d\sigma(x)\cdot\nu_j\leq B\|f\|^2
\end{equation}
for any $f\in\mathcal L^2(\mathcal S^n)$. For the wavelet transform we have
$$
\mathcal W_\Psi f(b_j,x)=\sum_{l=0}^\infty\overline{\gamma\bigl(b_j\tau(l)\bigr)}\,f_l(x),
$$
hence, by the Parseval identity~\eqref{eq:Parseval} we obtain
$$
\int_{\mathcal S^n}|\mathcal W_\Psi f(b_j,x)|^2\,d\sigma(x)=\sum_{l=0}^\infty\left|\gamma\bigl(b_j\tau(l)\bigr)\right|^2\sum_{\kappa=1}^{N(n,l)}|a_l^\kappa(f)|^2,
$$
and since all the sums converge absolutely, we may change the order of summation with respect to~$j$ and with respect to~$l$. Consequently, we may write~\eqref{eq:semifr_cond} as
\begin{equation}\label{eq:nier}\begin{split}
A\sum_{l=0}^\infty\sum_{\kappa=1}^{N(n,l)}|a_l^\kappa(f)|^2
   &\leq\sum_{l=0}^\infty\sum_{j=0}^\infty\left|\gamma\bigl(b_j\tau(l)\bigr)\right|^2\nu_j\sum_{\kappa=1}^{N(n,l)}|a_l^\kappa(f)|^2\\
&\leq B\sum_{l=0}^\infty\sum_{\kappa=1}^{N(n,l)}|a_l^\kappa(f)|^2.
\end{split}\end{equation}
Now, for any $l\in\mathbb N_0$ set $f=Y_l^{(0,0,\dots,0)}$, i.e., the zonal hyperspherical harmonic of degree~$l$. If $\kappa=1$ is the index of the sequence $(0,0,\dots,0)$, then we have
$$
\alpha_l^1(f)=1,
$$
and all the other coefficients~$\alpha_l^\kappa$ vanish. Thus, \eqref{eq:semiframe_cond}  follows from~\eqref{eq:nier}.

On the other hand, suppose~\eqref{eq:semiframe_cond} holds, then by the Parseval identity and Funk--Hecke formula we obtain from~\eqref{eq:nier} the inequality~\eqref{eq:semifr_cond}, and hence, the set~\eqref{eq:set_semifr} is a semi--frame.\hfill$\Box$

Note that this is in principle the same result as that obtained by Bogdanova \emph{et al.} in~\cite{BVAJM} and by Antoine in~\cite{jpA} (for wavelets defined in~\cite{AV} and not the ones we are working with; a relationship between the two wavelet definitions for the two--dimensional sphere is discussed in~\cite{ADJV02}, in~\cite{AMVA08}, and in~\cite{IIN14CWT}). In the case of wavelets derived from an approximate identity and satisfying~\eqref{eq:pseudo_generating}, we can simplify the condition~\eqref{eq:semiframe_cond}, as well as show that the frame~\eqref{eq:set_semifr} can be arbitrarily close to a tight one.

\begin{thm}\label{thm:semiframes}Let~$\{g_a:\,a\in\mathbb{R}_+\}$ be a wavelet family with \eqref{eq:pseudo_generating} for~$\gamma$ such that $\int_0^\infty\left||\gamma^2|^\prime(t)\right|\,dt<\infty$. Then, for any $\epsilon>0$ and $\delta>0$ there exists a constant~$\rho$ such that for any sequence~$\mathcal{B}=(b_j)_{j\in\mathbb{N}_0}$ with $b_0\geq-\log\rho$ and $1<b_{j}/b_{j+1}<1+\rho\cdot\delta$ the family \mbox{$\{g_{b_j,x},\,b_j\in\mathcal{B},\,x\in\mathcal S^n\}$} is a semi-continuous frame for~$\mathcal{L}^2(\mathcal S^n)$, satisfying the frame condition~\eqref{eq:framedef} with the prescribed~$\epsilon$.
\end{thm}

\begin{bfseries}Proof. \end{bfseries}Analogous to the proof of \cite[Theorem~5]{IH10}.\hfill$\Box$

It is easy to verify that Abel--Poisson wavelets with Gegenbauer coefficients of the kernel
$$
\widehat{\Psi_\rho}(l)=\frac{l+\lambda}{\lambda}\,\sqrt{2l\rho}\,e^{-l\rho},\qquad l\in\mathbb N_0,
$$
and Gauss--Weierstrass wavelets given by
$$
\widehat{\Psi_\rho}(l)=\frac{l+\lambda}{\lambda}\,\sqrt{2l(l+2\lambda)\rho}\,e^{-l(l+2\lambda)\rho},\qquad l\in\mathbb N_0,
$$
as well as Poisson multipole wavelets satisfy the conditions of the above theorem.

\section{Discrete frames of spherical wavelets}\label{sec:discrete_frames}

In this section we present a general theorem which states that under certain conditions on the localization of the reproducing kernel of the wavelet transform and its surface gradient fully discrete frames exist. It is worth noting that the phase--space discretization is on the one hand irregular, while on the other hand the density of sampling point distribution is quite uniform. It is an advantage of our approach compared to the equiangular grid used for the $2$--sphere in~\cite{jpA} and \cite{BVAJM}. Although equiangular grids allow the use of the Fourier transform, a concentration of points around the poles is their big drawback, such that a generalization to higher dimension is not reasonable.

The proof is based on the perturbation result from Corollary~\ref{cor:discretization}, and the existence of semi--continuous scale--discrete frames is assumed. The proof is analogous to that given in~\cite{IH10} for the two-dimensional case, but since a careful estimation of error estimating integrals is needed, we write here all the details.

The discretization is performed as in the following definition.

\begin{df}
We say a \emph{grid} $\Lambda \in \mathbb R_+\times\mathcal S^n $ is \emph{of a type} $(\delta,\Xi)$ \index{grid of type $(X,Y,\delta )$} if the following holds: There is a
sequence of scales $\mathcal{B}=(b_{j})_{j\in \mathbb{N}_{0}}$ such that the ratio $b_{j}/b_{j+1}$ is uniformly bounded from below and from above with
the lower bound larger than~$1$ and the upper bound equal to~$\delta$
$$
\widetilde{\delta}\leq b_{j}/b_{j+1}\leq \delta,\quad \widetilde{\delta}>1.
$$
At each scale~$b=b_{j}$, there is a measurable partition of~$\mathcal S^n$ $\mathcal{P}_{b}=\{\mathcal{O}_{k}^{(b)},\,k=1,2,\dots ,K_{b}\}$ into simply connected sets such that the diameter of each set (measured in geodesic distance) is not larger than~$\Xi\, b$. Each of these sets contains exactly one point of the grid.
\end{df}

\begin{thm}\label{thm:Pframes}Let $\mathcal{B}\subset \mathbb{R}_+$ be a set of scales $b_j$, $j\in\mathbb N_0$, with
\begin{equation*}
\widetilde\delta\leq b_j/b_{j+1}\leq\delta\quad\text{for some}\quad \widetilde\delta>1
\end{equation*}
 and $\nu(b_j)=\nu_j=\log ( b_j / b_{j+1})$. Further, suppose that $\{\Pi(a,x;\cdot,\cdot), (a,x)\in  \mathcal{B}\times\mathcal S^n\}$ is a frame with respect to the weights $\{C\nu_j\delta_{b_j}\}$ with a constant~$C$,
\begin{equation*}
(1-\Delta) \|s\|^2 \leq C\sum_j \int_{\mathcal S^n} \left| \int_{ \mathbb R_+\times\mathcal S^n }\hspace{-1em}
   \Pi(b_j,x; b,y)\, s(b,y) \,d\sigma(y) \frac{db}{b} \right|^2 d\sigma(x)\,\nu_j   \leq (1+\Delta) \|s\|^2.
\end{equation*}
 If in addition the reproducing kernel~$\Pi$ satisfies
\begin{equation}  \label{eq:extkern3}
\left.
\begin{array}{c}
|\Pi(a,x;b,y)| \\[0.5em]
\left|(a+b)\,\nabla_{\!\ast}\Pi(a,x;b,y)\right|
\end{array}
\right\} \leq(ab)^{n+\epsilon}\cdot
\begin{cases}
\frac{\mathfrak{d}}{(a+b)^{3n+2\epsilon}}, & \hspace{-0.7em}_{\angle(x,y)\leq\omega[a+(2-\tilde\epsilon)b],} \\
\frac{\mathfrak{d}}{\angle(x,y)^{3n+2\epsilon}}, & \hspace{-0.7em}_{\angle(x,y)>\omega(a+\tilde\epsilon b)},
\end{cases}
\end{equation}
for $a,b\leq b_0$ and for some positive constants~$\mathfrak{d}$, $\omega$, $\epsilon$ and $\tilde\epsilon<1/2$, where~$\nabla_\ast$~is the surface gradient with respect to any of the variables~$x$ or~$y$, then there exists a constant~$\rho$, such that for any grid~$\Lambda\subset \mathcal B\times\mathcal S^n $ of a type $(\delta,\Xi)$ with $\Xi\leq \rho$ the family $\{\Pi(b,y;\cdot,\cdot), (b,y)\in\Lambda\}$ is a frame with weight $C\sum\mu(b,y) \,\delta_b\, \delta_y$ for $\mu(b,y)=\sigma(\mathcal{O}_k^{(b)})$, $y\in\mathcal{O}_k^{(b)}$.
\end{thm}

The proof makes use of a convolution estimate for functions over the parameter space $\mathbb{H}= \mathbb R_+\times\mathcal S^n $. First we need a lemma, which is analogous to Young inequality for~$\mathbb{R}^n$. It is proven in the same way as \cite[Lemma~1]{IH10}.

\begin{lem}\label{lem:sphYoung}
Denote by~$\mathbb{K}$ the space~$\mathbb{R}_+\times\mathbb{R}_+$ with the measure $(da/a,\theta^{n-1}\,d\theta)$. Let~$F$ be such a
function $\mathbb{H}\times\mathbb{H}\to\mathbb{R}$ that
\begin{equation*}
F(a,x;b,y)=\frac{1}{b^n}\cdot f\left(\frac{a}{b},\frac{\angle(x,y)}{b}\right), \quad f\in\mathcal{L}^1(\mathbb{K}),
\end{equation*}
and $T\in\mathcal{L}^p(\mathbb{H})$, $p\geq1$. Then the following holds
\begin{equation*}
\|F\circ T\|_{\mathcal{L}^p(\mathbb{H})}\leq \Sigma_{n-1}\,\|f\|_{\mathcal{L}^1(\mathbb{K})}\cdot\|T\|_{\mathcal{L}^p(\mathbb{H})},
\end{equation*}
where the operation~$\circ$ is defined by
\begin{equation*}
F\circ T(a,x)=\int_{\mathbb{H}}\,F(a,x;b,y)\,T(b,y)\,\,d\sigma(y)\frac{db}{b}.
\end{equation*}
\end{lem}

\noindent\emph{Proof. } Let~$R$ be a non-negative function in~$\mathcal{L}^q(\mathbb{H})$ with $p^{-1}\!+q^{-1}\!=1$. We may also suppose, that~$F$ and~$T$ are non-negative. Then
\begin{align*}
\left<F\circ T,R\,\right>&=\int_{\mathbb{H}} F\circ T(a,x)\,R(a,x)\,\,d\sigma(x)\frac{da}{a} \\
&=\int_{\mathbb{H}}\int_{\mathbb{H}}\frac{1}{b^n}f\!\left(\frac{a}{b},\frac{\angle(x,y)}{b}\right)T(b,y)\,R(a,x)\,\,d\sigma(y)\frac{db}{b}\,\,d\sigma(x)\frac{da}{a}.
\end{align*}
By the change of the variables $a/b\mapsto a$ and exchanging the integrals (since all functions are positive, the integrals may only converge absolutely) we obtain
\begin{equation*}
\left<F\circ T,R\,\right>=\int_{\mathbb{H}}\int_{\mathbb{H}}\frac{1}{b^n}f\!\left(a,\frac{\angle(x,y)}{b}\right)R(a,xb)\,\,d\sigma(x)\frac{da}{a}\,T(b,y)\,\,d\sigma(y)\frac{db}{b}.
\end{equation*}
Consider the inner integral with respect to $d\sigma(x)$, which for simplicity purposes we write as $\int g(x\cdot y)\,r(x)\,d\sigma(x)$. Let $A=A_y$ be an isometry of the sphere which maps~$y$ to the North Pole~$\hat{e}$ and~$\hat{e}$ to~$y$. Then
\begin{equation}\label{eq:integral_isometry}
\int g(x\cdot y)\,r(x)\,d\sigma(x)=\int g(Ax\cdot y)\,r(Ax)\,d\sigma(Ax)=\int g(x\cdot A^\ast y)\,r(Ax)\,d\sigma(x).
\end{equation}
Since
$$
\hat e\cdot A^\ast y=Ay\cdot\hat e=\hat e\cdot\hat e=1,
$$
we have $A^\ast y=\hat e$ and~\eqref{eq:integral_isometry} yields
$$
\int g(x\cdot y)\,r(x)\,d\sigma(x)=\int g(x\cdot\hat{e})\,r(Ax)\,d\sigma(x).
$$
Now, $Ax$ describes the position of the point~$x$ relative to the point~$y$ (depending also on the position of the North Pole). Let~$x$ be fixed; by~$\mathcal{R}_x $ we denote the function \mbox{$(y,a)\mapsto R(A_ya,x)\left(=r(Ax)\right)$}. Since~$A$ was an isometry, we have
\begin{equation*}
\int_{\mathcal S^n} \mathcal{R}_x(y,a)\,d\sigma(y)=\int_{\mathcal S^n} R(y,a)\,d\sigma(y).
\end{equation*}
Then we have (once again exchanging the integrals)
\begin{equation*}
\left<F\circ T,R\,\right>=\int_{\mathbb{H}}\int_{\mathbb{R}_+}\!\int_{\mathcal S^n} \mathcal{R}_x(y,ab)\,T(b,y)\,d\sigma(y)\,
   \frac{1}{b^n}f\!\left(a,\frac{\theta}{b}\right)\frac{db}{b}\,\,d\sigma(x)\frac{da}{a},
\end{equation*}
where $\theta=\angle(x,\hat{e})$, and further, by H\"older inequality,
\begin{equation*}
\left<F\circ T,R\,\right>\leq\int_{\mathbb{H}}\int_{\mathbb{R}_+}\|R(\cdot,ab)\|_{\mathcal{L}^q(\mathcal S^n)}\|T(\cdot,b)\|_{\mathcal{L}^p(\mathcal S^n)}
   \frac{1}{b^n}f\!\left(a,\frac{\theta}{b}\right)\frac{db}{b}\,\,d\sigma(x)\frac{da}{a}.
\end{equation*}
Now, the integral over~$\mathcal S^n$ may be estimated as follows:
\begin{align*}
\int_{\mathcal S^n}\frac{1}{b^n}&f(a,\theta/b)\,d\sigma(x)=\Sigma_{n-1}\int_0^\pi\frac{1}{b^n}f(a,\theta/b)\sin^{n-1}\theta\,d\theta\\
&=\Sigma_{n-1}\int_0^{\pi/b}\!f(a,\theta)\,\frac{\sin^{n-1}(b\theta)}{b^{n-1}}\,d\theta\leq\Sigma_{n-1}\int_0^{\pi/b}\!f(a,\theta)\,\theta^{n-1}\,d\theta\\
&\leq\Sigma_{n-1}\int_0^\infty\!f(a,\theta)\,\theta^{n-1}\,d\theta=\Sigma_{n-1}\|f(a,\cdot)\|_{\mathcal{L}^1(\mathbb{R}_+,\theta^{n-1} d\theta)},
\end{align*}
and therefore, by H\"older inequality with respect to $db/b$,
\begin{align*}
\left<F\circ T,R\,\right>&\leq\Sigma_{n-1}\int_{\mathbb{R}_+}\!\int_{\mathbb{R}_+}\!\|f(a,\cdot)\|_{\mathcal{L}^1(\mathbb{R}_+,\theta^{n-1}d\theta)}
   \|T(\cdot,b)\|_{\mathcal{L}^p(\mathcal S^n)} \|R(\cdot,ab)\|_{\mathcal{L}^q(\mathcal S^n)}\frac{da}{a}\frac{db}{b} \\
&\leq\Sigma_{n-1}\int_{\mathbb{R}_+}\!\|f(a,\cdot)\|_{\mathcal{L}^1(\mathbb{R}_+,\theta^{n-1}d\theta)}\frac{da}{a}
   \cdot\|T\|_{\mathcal{L}^p(\mathbb{H})}\|R\|_{\mathcal{L}^q(\mathbb{H})} \\
&=\Sigma_{n-1}\,\|f\|_{\mathcal{L}^1(\mathbb{K})}\|T\|_{\mathcal{L}^p(\mathbb{H})}\|R\|_{\mathcal{L}^q(\mathbb{H})}.
\end{align*}
Therefore, we have by Riesz representation theorem
\begin{equation*}
\|F\circ T\|_{\mathcal{L}^p(\mathbb{H})}\leq \Sigma_{n-1}\,\|f\|_{\mathcal{L}^1(\mathbb{K})}\|T\|_{\mathcal{L}^p(\mathbb{H})}.
\end{equation*}
Since by assumption all the norms are finite, the exchanges of integrals were justified. \hfill $\Box$

\noindent\emph{Proof of Theorem~\ref{thm:Pframes}. } According to Lemma~\ref{lem:sphYoung} and Corollary~\ref{cor:discretization} it is enough to show that
\begin{equation*}
\begin{split}
D&=\left|\sum_{(b,y)\in\Lambda}\Pi(a,x;b,y)\,\Pi(b,y;c,z)\,\mu(b,y)\right. \\
&\quad\left.-C\sum_{b\in\mathcal B}\int_{{\mathcal S^n}}\Pi(a,x;b,y)\,\Pi(b,y;c,z)\,d\sigma(y)\,\nu(b)\right|
\end{split}
\end{equation*}
is less than
$$
\frac{1}{C}\cdot\widetilde\Delta\cdot\frac{1}{c^n}\,f\left(\frac{a}{c},\frac{\angle(x,z)}{c}\right)
$$
for some $f\in\mathcal{L}^1(\mathbb{K})$ with $\|f\|=\frac{1}{\Sigma_{n-1}}$ and $\widetilde\Delta\in(0,1-\Delta)$.

For fixed $(a,x)$, $(c,z)$ and~$b\in\mathcal B$, set $F(y)=\Pi(a,x;b,y)$ and $G(y)=\Pi(c,z;b,y)$. Let~$\mathcal{K}^x$ denote the set of points, where $F$ is 'large', i.e., $\mathcal{K}^x=\{y\in{\mathcal S^n}:\,\angle(x,y)\leq\omega(a+b)\} $. Similarly, denote by~$\mathcal{K}^z$ the set '$G$ large', i.e., $\mathcal{K}^z=\{y\in{\mathcal S^n}:\,\angle(y,z)\leq\omega(c+b)\}$. If the sets~$\mathcal{K}^x$ and~$\mathcal{K}^z$ are not disjoint, we split the error that one makes by exchanging integration over~${\mathcal S^n}$ by summation over $\{y\in{\mathcal S^n}:\,(b,y)\in\Lambda\}$ into two parts

\begin{itemize}
\item[--] $I_1(b)$: $F(y)$ 'large' or $G(y)$ 'large', i.e., over the set $\mathcal{D}=\mathcal{K}^x\cup\mathcal{K}^z$;

\item[--] $I_4(b)$: $F(y)$ 'small' and $G(y)$ 'small', i.e., for $\mathcal{G}={\mathcal S^n}\backslash(\mathcal{K}^x\cup\mathcal{K}^z)$.
\end{itemize}

In the other case, if the sets~$\mathcal{K}^x$ and~$\mathcal{K}^z$ have an empty intersection, we consider three parts:

\begin{itemize}
\item[--] $I_2(b)$: $F(y)$ 'large', $G(y)$ 'small', i.e., for $\mathcal{E}=\mathcal{K}^x$;

\item[--] $I_3(b)$: $F(y)$ 'small', $G(y)$ 'large', i.e., for $\mathcal{F}=\mathcal{K}^z$;

\item[--] $I_4(b)$: $F(y)$ 'small', $G(y)$ 'small', i.e., for $\mathcal{G}={\mathcal S^n}\backslash(\mathcal{K}^x\cup\mathcal{K}^z)$.
\end{itemize}

Each of the errors may be estimated in the following way: for every set~$\mathcal{O}=\mathcal{O}_k^{(b)}$ the difference between the highest and the lowest value of~$F(\eta)\cdot G(\eta)$, $\eta\in\mathcal{O}$, is less than or equal to
\begin{equation*}
\sup_{\eta\in\mathcal{O}}|\nabla_{\!\!\ast}\left[F(\eta)\cdot G(\eta)\right]|\cdot\text{diam}(\mathcal{O}),
\end{equation*}
and hence the difference between $\int_{\mathcal{O}}F(\eta)\,G(\eta)\,d\sigma(\eta)\,\nu(b)$ and $F(y)\,G(y)\,\mu(b,y)$ for $y=y_k^{(b)}$ is less than or equal to
\begin{equation}  \label{eq:esterrwin}
\left(\sup_{\eta\in\mathcal{O}}|\nabla_{\!\!\ast}F(\eta)|\cdot\sup_{\eta\in\mathcal{O}}|G(\eta)|
   +\sup_{\eta\in\mathcal{O}}|F(\eta)|\cdot\sup_{\eta\in\mathcal{O}}|\nabla_{\!\!\ast}G(\eta)|\right) \cdot\text{diam}(\mathcal{O})\cdot\mu(b,y).
\end{equation}

When summing up over all the sets~$\mathcal{O}$ that have a non-empty intersection with one of the sets ($\mathcal{D}$, $\mathcal{E}$, $\mathcal{F}$ or $\mathcal{G}$), we may calculate suprema over the whole set and choose the largest possible~$\text{diam}(\mathcal{O})$ ($=\Xi b$), which we denote by~$R$. The sum of~$\mu(b,y)$ is then not larger than the volume of~$\mathcal{Y}_{R}$ multiplied by~$\nu(b)$, with~$\mathcal{Y}_{R}$ denoting the $R$--parallel extension of~$\mathcal{Y}$, i.e.,
\begin{equation*}
\mathcal{Y}_{R}=\{\left.\eta\in{\mathcal S^n}\right|\,\exists y\in\mathcal{Y}:\,\angle(\eta,y)\leq R\},
\end{equation*}
where~$\mathcal{Y}$ means one of the sets~$\mathcal{D}$, $\mathcal{E}$, $\mathcal{F}$ or $\mathcal{G}$.

We introduce the notation $\alpha=a/c$, $\beta=b/c$, $\theta=\angle(x,z)$, $\vartheta=\theta/c$ and $f_\iota(\alpha,\vartheta)=\sum_b c^n I_\iota(b)$ for $\iota=1,2,3,4$ and $b\in\mathcal{B}$, but possibly not all the scales. The constant~$\mathfrak{c}$ may change its value from line to line.

{\bfseries Part 1) } For~$I_1$ we have
\begin{equation}  \label{eq:est}
c^nI_1(b)\leq c^n\cdot\mathfrak{c}\cdot\left(\frac{1}{a+b}+\frac{1}{c+b}\right) \cdot\frac{(ab)^{n+\epsilon}}{(a+b)^{3n+2\epsilon}}
   \cdot\frac{(bc)^{n+\epsilon}}{(c+b)^{3n+2\epsilon}}\cdot \Xi b\cdot\sigma(\mathcal{D}_{R})\cdot\nu(b).
\end{equation}
The set~$\mathcal{D}$ is contained in~$\mathcal{K}^x\cup\mathcal{K}^z$ and hence, the volume of~$\mathcal{D}_{R}$ is bounded by \mbox{$2\cdot\text{vol}(\mathcal{K}_{R})$}, where~$\mathcal{K}$ is the larger of the balls~$\mathcal{K}^x$ and~$\mathcal{K}^z$. This is given by $\mathfrak c\cdot\omega^n\cdot\left(c+(1+\Xi/\omega)\,b\right)^n\leq\mathfrak{c}\,(c+b)^n$ if $a<c$, respectively $\mathfrak c\cdot\omega^n\cdot\left(a+(1+\Xi/\omega)\,b\right)^n\leq\mathfrak{c}(a+b)^n$ if $a\geq c$. If $\alpha\leq1$, we obtain from~\eqref{eq:est}:
\begin{equation}  \label{eq:est1asmall}
c^nI_1(b)\leq\mathfrak{c}\,\Xi\cdot\frac{\alpha^{n+\epsilon}}{(\alpha+\beta)^{n+\epsilon/2}} \cdot\frac{\beta^{2n+1+3\epsilon/2}}{(\alpha+\beta)^{2n+1+3\epsilon/2}}    \cdot\frac{\beta^{\epsilon/2}}{(1+\beta)^{2n+2\epsilon}}\cdot\nu(b).
\end{equation}
The second fraction is smaller than~$1$, and the last one ensures the summability over~$b$, thus, we have the estimation
\begin{equation}
f_1\leq\mathfrak{c}\,\Xi\cdot\alpha^{\epsilon/2},  \tag{\cAa}
\end{equation}
which we use for the error estimation for $\vartheta\leq\omega(\alpha+1)$. For large~$\vartheta$, $\vartheta>\omega(\alpha+1)$, we need a sharper result. Since the sets~$\mathcal{K}^x$ and~$\mathcal{K}^z$ have a non--empty intersection only for~$b$ such that $\omega(\alpha+2\beta+1)\geq\vartheta$, i.e., $2(1+\beta)\geq\vartheta/\omega+1-\alpha$, we may enlarge the last fraction in the estimation~\eqref{eq:est1asmall}, and write
\begin{equation*}
\frac{\beta^{\epsilon/2}}{(1+\beta)^{2n+2\epsilon}}\leq\frac{\beta^{\epsilon/2}}{(1+\beta)^{2\epsilon}} \cdot\frac{\mathfrak{c}}{[\vartheta+\omega(1-\alpha)]^{2n}}.
\end{equation*}
Consequently, we obtain
\begin{equation}
f_1\leq\mathfrak{c}\,\Xi\cdot\frac{\alpha^{\epsilon/2}}{[\vartheta+\omega(1-\alpha)]^{2n}}.  \tag{\cBa}
\end{equation}

In the other case, $\alpha>1$, we get
\begin{equation*}
c^nI_1(b)\leq\mathfrak{c}\,\Xi\cdot\frac{\alpha^{n+\epsilon}}{(\alpha+\beta)^{2n+2\epsilon}} \cdot\frac{\beta^{2n+1+2\epsilon}}{(1+\beta)^{3n+1+2\epsilon}}\cdot\nu(b).
\end{equation*}
For $\vartheta\leq\omega(1+\alpha)$ we then have
\begin{equation}
f_1\leq\mathfrak{c}\,\Xi\cdot\frac{1}{\alpha^{n+\epsilon}}  \tag{\cCa}
\end{equation}
and for $\vartheta>\omega(1+\alpha)$ we write
\begin{equation}
f_1\leq\mathfrak{c}\,\Xi\cdot\frac{1}{(\alpha+\beta)^{n+\epsilon}}\leq\mathfrak{c}\Xi \cdot\frac{1}{\alpha^{\epsilon/2}\,[\vartheta+\omega(\alpha-1)]^{n+\epsilon/2}},  \tag{\cDa}
\end{equation}
since for~$b$ we take into account the relation $2(\alpha+\beta)\geq\vartheta/\omega+\alpha-1$.

{\bfseries Part 2) } In the second case, $I_2(b)$, we consider only the scales for which $\mathcal{K}^x$ and $\mathcal{K}^z$ have an empty intersection, i.e., $b$~such that \mbox{$\vartheta>\omega(\alpha+2\beta+1)$}. For the error made in the whole set~$\mathcal{E}$ we use the formula~\eqref{eq:esterrwin} with~$\mu(b,y)$ replaced by the volume of~$\mathcal{E}_{R}$ (i.e., the volume of~$(\mathcal{K}^x)_{R}$) multiplied by~$\nu(b)$. The supremum of the modules of~$G$ and~$\nabla_{\!\!\ast}G$ is estimated by their values
at the point nearest~$\mathcal{K}^z$. Since we have to consider all the sets~$\mathcal{O}_k^{(b)}$ that have a non--empty intersection with~$\mathcal{E}$, we choose the angular argument in~\eqref{eq:extkern3} to be equal to~$\theta-\omega(a+b)-R$.

We have to assume that the maximum diameter of a partition set is less than $\mathfrak{c}\cdot\omega b$, with some~$\mathfrak{c}<1/2$. For the sake of simplicity, we set $R\leq\omega b/3$. Altogether we obtain
\begin{equation}  \label{eq:estI2}
\begin{split}
c^nI_2(b)&\leq c^n\cdot\mathfrak{c}\cdot\left(\frac{1}{a+b}+\frac{1}{c+b}\right) \\
&\cdot\frac{(ab)^{n+\epsilon}}{(a+b)^{3n+2\epsilon}} \cdot\frac{(bc)^{n+\epsilon}}{[\theta-\omega(a+4b/3)]^{3n+2\epsilon}}\cdot\Xi b\cdot(a+b)^n\cdot\nu(b).
\end{split}
\end{equation}
Further, in the considered range of scales we have $\vartheta/\omega>\alpha+2\beta+1$, and this inequality implies $\vartheta-\omega(\alpha+4\beta/3)>[\vartheta+\omega(2-\alpha)]/3$ as well as $\vartheta-\omega(\alpha+4\beta/3)>\mathfrak c\,(1+\beta)$.

For $\alpha\leq1$, we write the estimation~\eqref{eq:estI2} in the form
\begin{equation*}
c^nI_2(b)\leq\mathfrak{c}\,\Xi\cdot\frac{\alpha^{n+\epsilon}}{(\alpha+\beta)^{1+\epsilon}} \cdot\frac{\beta^{2n+\epsilon}}{(\alpha+\beta)^{2n+\epsilon}}
   \cdot\frac{\beta^{1+\epsilon}}{(1+\beta)^{1+2\epsilon}}\cdot\frac{1}{[\vartheta+\omega(2-\alpha)]^{3n-1}}\cdot\nu(b),
\end{equation*}
that yields
\begin{equation}
f_2\leq\mathfrak{c}\,\Xi\cdot\frac{\alpha^{n-1}}{[\vartheta+\omega(2-\alpha)]^{3n-1}},
\tag{\cAb}
\end{equation}
and for $\alpha>1$ we have
\begin{equation*}
c^nI_2(b)\leq\mathfrak{c}\,\Xi\cdot\frac{\alpha^{n+\epsilon}}{(\alpha+\beta)^{n+1+\epsilon}} \cdot\frac{\beta^{n-1+\epsilon}}{(\alpha+\beta)^{n-1+\epsilon}}
   \cdot\frac{\beta^{n+2+\epsilon}}{(1+\beta)^{n+2+2\epsilon}} \cdot\frac{1}{[\vartheta+\omega(2-\alpha)]^{2n-1}}\cdot\nu(b);
\end{equation*}
consequently,
\begin{equation}
f_2\leq\frac{\mathfrak{c}\,\Xi}{\alpha\,[\vartheta+\omega(2-\alpha)]^{2n-1}}.
\tag{\cBb}
\end{equation}

{\bfseries Part 3) } Similarly as in the previous case, we obtain from
\begin{equation*}\begin{split}
c^nI_3(b)&\leq c^n\cdot\mathfrak{c}\cdot\left(\frac{1}{a+b}+\frac{1}{c+b}\right) \\
&\cdot\frac{(ab)^{n+\epsilon}}{[\theta-\omega(c+4b/3)]^{3n+2\epsilon}} \cdot\frac{(bc)^{n+\epsilon}}{(c+b)^{3n+2\epsilon}}\cdot \Xi b\cdot(c+b)^n\cdot\nu(b)
\end{split}\end{equation*}
the estimations
\begin{equation*}
c^nI_3(b)\leq\mathfrak{c}\,\Xi\cdot\alpha^{n+\epsilon}\cdot\frac{1}{[\vartheta+\omega(2\alpha-1)]^{n+1}}
   \cdot\frac{\beta^{2n+2\epsilon}}{(\alpha+\beta)^{2n+2\epsilon}}\cdot\frac{\beta}{(1+\beta)^{2n+2\epsilon}}\cdot\nu(b),
\end{equation*}
for $\alpha\leq1$ and
\begin{equation*}
c^nI_3(b)\leq\mathfrak{c}\,\Xi\cdot\alpha^{n+\epsilon}\cdot\frac{1}{[\vartheta+\omega(2\alpha-1)]^{3n}}
   \cdot\frac{\beta^{2\epsilon}}{(\alpha+\beta)^{2\epsilon}}\cdot\frac{\beta^{2n+1}}{(1+\beta)^{2n+1+2\epsilon}}\cdot\nu(b)
\end{equation*}
for $\alpha>1$. They yield
\begin{equation}
f_2\leq\mathfrak{c}\,\Xi\cdot\frac{\alpha^{n+\epsilon}}{[\vartheta+\omega(2\alpha-1)]^{n+1}}  \tag{\cAc}
\end{equation}
for $\alpha\leq1$ and
\begin{equation}
f_2\leq\mathfrak{c}\,\Xi\cdot\frac{\alpha^{n+\epsilon}}{[\vartheta+\omega(2\alpha-1)]^{3n}}  \tag{\cBc}
\end{equation}
for $\alpha>1$.

{\bfseries Part 4) } a) Consider large~$\theta$ and small scales~$b$, that is, those satisfying the condition $\theta>\omega(a+2b+c)$. For the points~$y $ on the sphere that lie closer to the spherical ball $\mathcal{K}^x$, i.e., elements of the set
$$
\mathcal{R}_x:=\{y\in{\mathcal S^n}\backslash\mathcal{K}^x:\,\angle(x,y)-\omega a\leq\angle(z,y)-\omega c\},
$$
and for one set~$\mathcal{O}=\mathcal{O}_k^{(b)}$, we estimate the error using formula~\eqref{eq:esterrwin}; the terms $\sup_{\eta\in\mathcal{O}}|G(\eta)|$ and~$\sup_{\eta\in\mathcal{O}}|\nabla_{\!\!\ast}G(\eta)|$ may be replaced by the largest possible value in the $R$--parallel extension of~$\mathcal{R}_x$, i.e.
\begin{equation}  \label{eq:biggestg}
\sup_{\eta\in\mathcal{O}}|G(\eta)|\leq\frac{\mathfrak{c}\cdot(cb)^{n+\epsilon}}{\theta_z^{3n+2\epsilon}}\quad\text{resp.}\quad
   \sup_{\eta\in\mathcal{O}}|\nabla_{\!\!\ast}G(\eta)|\leq\frac{\mathfrak{c}\cdot(cb)^{n+\epsilon}}{(c+b)\cdot\theta_z^{3n+2\epsilon}}
\end{equation}
with
\begin{equation*}
\theta_z=\omega c+\frac{\theta-\omega(a+c)}{2}-R\geq\frac{\theta+\omega(2c-a)}{3}\geq\omega\left(c+\frac{2}{3}\,b\right).
\end{equation*}
Further, $\sup_{\eta\in\mathcal{O}}|\nabla_{\!\!\ast}F(\eta)|\cdot\mu(b,y)$ resp. $\sup_{\eta\in\mathcal{O}}|F(\eta)|\cdot\mu(b,y)$ may be estimated by
\begin{equation*}
\frac{(ab)^{n+\epsilon}}{a+b}\int_\mathcal{O}\frac{\,d\sigma(y)}{\left(\angle(x,y)-R\right)^{3n+2\epsilon}}\quad\text{resp.}
\quad(ab)^{n+\epsilon}\int_\mathcal{O}\frac{\,d\sigma(y)}{\left(\angle(x,y)-R\right)^{3n+2\epsilon}}
\end{equation*}
multiplied by~$\nu(b)$. The bound we obtain for the error is larger if we sum up over \emph{all} the partition sets having a non--empty intersection with the complement of~$\mathcal{K}^x$ (with $\sup_\eta|G(\eta)|$ given by~\eqref{eq:biggestg}, a property that does not hold in the whole~$({\mathcal S^n}\backslash\mathcal{K}^x)_{R}$). Since $R\leq\omega b/3$, we obtain
\begin{align*}
c^nI_4^{(x)}(b)&\leq c^n\cdot\mathfrak{c}\cdot\left(\frac{1}{a+b}+\frac{1}{c+b}\right) \\
&\cdot\int_{{\mathcal S^n}_x}\frac{(ab)^{n+\epsilon}}{\left(\angle(x,y)-\omega\,b/3\right)^{3n+2\epsilon}}
   \cdot\frac{(bc)^{n+\epsilon}}{\theta_z^{3n+2\epsilon}}\cdot \Xi b\,d\sigma(y)\cdot\nu(b)
\end{align*}
where~$I_4^{(x)}(b)$ means the error made in the set~$\mathcal{R}_x$ and ${\mathcal S^n}_x$ is the set $\{y\in{\mathcal S^n}:\,\angle(x,y)\geq\omega(a+2b/3)\}$.
Denote $\angle(x,y)$ by~$\chi$, then the integral is given by
$$
\int_{\omega(a+2b/3)}^\pi\frac{\Xi\,a^{n+\epsilon}\,b^{2n+1+2\epsilon}\,c^{n+\epsilon}} {(\chi-\omega\,b/3)^{3n+2\epsilon}\,\theta_z^{3n+2\epsilon}}\,
   \sin^{n-1}\chi\,d\chi.
$$
When replacing ~$\sin\chi$ by~$\chi=(\chi-\omega\,b/3)+\omega\,b/3 $ and the upper integration bound~$\pi$ by~$\infty$, we obtain
\begin{align*}
\int_{\omega(a+2b/3)}^\pi&\frac{\sin^{n-1}\chi\,d\chi}{(\chi-\omega\,b/3)^{3n+2\epsilon}}
   \leq\int_{\omega(a+b/3)}^\infty\frac{(\chi+\omega\,b/3)^{n-1}\,d\chi}{\chi^{3n+2\epsilon}}\\
&=\sum_{k=0}^{n-1}\binom{n-1}{k}\frac{(\omega\,b/3)^{n-1-k}}{(3n+2\epsilon-k-1)\left(\omega(a+b/3)\right)^{3n+2\epsilon-k-1}}\\
&\leq\frac{\mathfrak c}{(a+b/3)^{2n+2\epsilon}},
\end{align*}
and consequently
\begin{equation}\label{eq:estI4x}
c^nI_4^{(x)}(b)\leq\mathfrak{c}\,\Xi\cdot\left(\frac{1}{a+b}+\frac{1}{c+b}\right)
   \cdot\frac{a^{n+\epsilon}\,b^{2n+1+2\epsilon}\,c^{2n+\epsilon}} {(a+b/3)^{2n+2\epsilon}\,\theta_z^{3n+2\epsilon}}\cdot\nu(b)
\end{equation}
For $\alpha\leq1$ we can write:
\begin{equation}  \label{eq:estA}
c^nI_4^{(x)}(b)\leq\mathfrak{c}\,\Xi\cdot\frac{\alpha^n}{\alpha+\beta} \cdot\frac{\alpha^\epsilon\,\beta^{2n+\epsilon}}{(\alpha+\beta/3)^{2n+2\epsilon}}
   \cdot\frac{\beta^{1+\epsilon}}{\vartheta_z^{1+2\epsilon}}\cdot\frac{1}{\vartheta_z^{3n-1}}\cdot\nu(b).
\end{equation}
In the second case, $\alpha>1$, the inequality~\eqref{eq:estI4x} yields
\begin{equation}  \label{eq:estB}
c^nI_4^{(x)}(b)\leq\mathfrak{c}\,\Xi\cdot\frac{\alpha^{n+\epsilon}}{(\alpha+\beta/3)^{n+1+\epsilon}} \cdot\frac{\beta^{n+\epsilon}}{(1+\beta)\,(\alpha+\beta/3)^{n-1+\epsilon}}
   \cdot\frac{\beta^{n+1+\epsilon}}{\vartheta_z^{n+1+2\epsilon}}\cdot\frac{1}{\vartheta_z^{2n-1}}\cdot\nu(b).
\end{equation}

Analogously, for points closer to the other spherical ball, i.e., elements of
$$
\mathcal{R}_z:=\{y\in{\mathcal S^n}\backslash\mathcal{K}^z:\,\angle(x,y)-\omega a>\angle(z,y)-\omega c\},
$$
we obtain
\begin{equation}  \label{eq:estI4zint}
\begin{split}
c^nI_4^{(z)}(b)&\leq\mathfrak{c}\cdot\left(\frac{1}{a+b}+\frac{1}{c+b}\right)\\
&\cdot\int_{{\mathcal S^n}_z}\frac{(ab)^{n+\epsilon}}{\theta_x^{3n+2\epsilon}} \cdot\frac{(bc)^{n+\epsilon}}{\left(\angle(z,y)-R\right)^{3n+2\epsilon}}
   \cdot c^n\cdot \Xi b\,d\sigma(y)\cdot\nu(b),
\end{split}
\end{equation}
where ${\mathcal S^n}_z=\{y\in{\mathcal S^n}:\,\angle(z,y)\geq\omega(c+2b/3)\}$ and
\begin{equation*}
\theta_x=\omega a+\frac{\theta-\omega(a+c)}{2}-R\geq\frac{\theta+\omega(2a-c)}{3}\geq\omega\left(a+\frac{2}{3}\,b\right)
\end{equation*}
(and~$I_4^{(z)}$ is the error made in the set~$\mathcal{R}_z$). The right--hand--side of the inequality~\eqref{eq:estI4zint} may be enlarged so that we get
$$
c^nI_4^{(z)}(b)\leq\mathfrak{c}\,\Xi\cdot\left(\frac{1}{\alpha+\beta}+\frac{1}{1+\beta}\right)
   \cdot\frac{\alpha^{n+\epsilon}\beta^{2n+1+2\epsilon}}{\vartheta_x^{3n+2\epsilon}\,(1+\beta/3)^{2n+2\epsilon}}\cdot\nu(b),
$$
and we write it for $\alpha\leq1$ as
\begin{equation}  \label{eq:estC}
c^nI_4^{(z)}(b)\leq\mathfrak{c}\,\Xi\cdot\frac{\alpha^{n+\epsilon}}{\vartheta_x^{n+1}} \cdot\frac{\beta^{2n+2\epsilon}}{(\alpha+\beta)\,\vartheta_x^{2n-1+2\epsilon}}
   \cdot\frac{\beta}{(1+\beta/3)^{2n+2\epsilon}}\cdot\nu(b).
\end{equation}
If $\alpha>1$, we use the factorization
\begin{equation}  \label{eq:estD}
c^nI_4^{(z)}(b)\leq\mathfrak{c}\,\Xi\cdot\frac{\alpha^{n+\epsilon}}{\vartheta_x^{3n}} \cdot\frac{\beta^{1+2\epsilon}}{(1+\beta)\,\vartheta_x^{2\epsilon}}
   \cdot\frac{\beta^{2n}}{(1+\beta/3)^{2n+2\epsilon}}\cdot\nu(b).
\end{equation}
b) If $\theta>\omega(a+c)$ and~$b$ is such that \mbox{$\theta\leq\omega(a+2b+c)$}, we estimate the error in a similar way, but we set
\begin{equation}  \label{eq:thxz}
\theta_x=\omega(a+b)-R\quad\text{and}\quad\theta_z=\omega(c+b)-R.
\end{equation}
We obtain again the estimations~\eqref{eq:estA}, \eqref{eq:estB}, \eqref{eq:estC} and \eqref{eq:estD}. In the first two of them, the denominator of the third fraction is always larger than or equal to powered \mbox{$\omega(1+2\beta/3)$}, and hence it ensures the summability over~$b$; the second fraction is not larger than a constant. In the inequalities~\eqref{eq:estC} and~\eqref{eq:estD}, one can replace the second fraction by a constant, since $\vartheta_x\geq\omega(\alpha+2\beta/3)$. Further, the estimations
\begin{equation*}
\theta_x\geq\frac{\theta+\omega(2a-c)}{3}\quad\text{and}\quad\theta_z\geq\frac{\theta+\omega(2c-a)}{3}
\end{equation*}
are also valid for~$\theta_x$ and~$\theta_z$ defined by~\eqref{eq:thxz} if the range of scales is bounded by $\theta/\omega\leq a+2b+c$, which is the case here. Consequently, we obtain from~\eqref{eq:estA} and~\eqref{eq:estC}
\begin{equation}
f_4\leq\mathfrak{c}\,\Xi\cdot\frac{\alpha^{n-1}}{[\vartheta+\omega(2-\alpha)]^{3n-1}} +\mathfrak{c}\cdot\frac{\alpha^{n+\epsilon}}{[\vartheta+\omega(2\alpha-1)]^{n+1}}
   \tag{\cAd}
\end{equation}
for $\alpha\leq1$ and from \eqref{eq:estB} and \eqref{eq:estD}
\begin{equation}
f_4\leq\mathfrak{c}\,\Xi\cdot\frac{1}{\alpha\,[\vartheta+\omega(2-\alpha)]^{2n-1}} +\mathfrak{c}\cdot\frac{\alpha^{n+\epsilon}}{[\vartheta+\omega(2\alpha-1)]^{3n}}\tag{\cBd}
\end{equation}
for $\alpha>1$.

c) Now, for $\theta\leq\omega(a+c)$, the sets~$\mathcal{K}^x$ and~$\mathcal{K}^z$ have a non--empty intersection for all scales~$b$. Since ${\mathcal S^n}\backslash(\mathcal{K}^x\cup\mathcal{K}^z)\subseteq{\mathcal S^n}\backslash\mathcal{K}^z$ and $\sup_{\eta\in({\mathcal S^n}\backslash\mathcal{K}
^z)_{R}}|G(\eta)|=|G(\omega(c+b)-R)|$, the inequality~\eqref{eq:estI4x} with $\theta_z\geq\omega(c+2b/3)\geq\mathfrak{c}\,(c+b)$:
\begin{equation*}
c^nI_4(b)\leq\mathfrak{c}\,\Xi\cdot\left(\frac{1}{\alpha+\beta}+\frac{1}{1+\beta}\right)
   \cdot\frac{\alpha^{n+\epsilon}\beta^{2n+1+2\epsilon}} {(\alpha+\beta/3)^{2n+2\epsilon}(1+\beta)^{3n+2\epsilon}}\cdot\nu(b)
\end{equation*}
yields an estimation of the error made in the whole set~$I_4$. For $\alpha\leq1$ we write it as
\begin{equation*}
c^nI_4(b)\leq\mathfrak{c}\,\Xi\cdot\frac{\alpha^{n+\epsilon}}{\alpha+\beta} \cdot\frac{\beta^{2n+2\epsilon}}{(\alpha+\beta/3)^{2n+2\epsilon}}
   \cdot\frac{\beta}{(1+\beta)^{3n+2\epsilon}}\cdot\nu(b)
\end{equation*}
and obtain for the sum over all scales:
\begin{equation}
f_4\leq\mathfrak{c}\,\Xi\cdot\alpha^{n-1+\epsilon}.  \tag{\cCd}
\end{equation}
In the opposite case, $\alpha>1$, one has
\begin{equation*}
c^nI_4(b)\leq\mathfrak{c}\,\Xi\cdot\frac{\alpha^{n+\epsilon}}{(\alpha+\beta/3)^{2n+2\epsilon}}\cdot\frac{\beta^{2n+1+2\epsilon}}{(1+\beta)^{3n+1+2\epsilon}}\cdot\nu(b),
\end{equation*}
and consequently
\begin{equation}
f_4\leq\mathfrak{c}\,\Xi\cdot\frac{1}{\alpha^{n+\epsilon}}.  \tag{\cDd}
\end{equation}
\newline
The following table sorts the obtained estimations:\vspace{0.5em}

\begin{center}
\begin{tabular}{|r|cc|}
\hline
\rule[-1em]{0em}{2.5em} & $\vartheta\leq\omega(\alpha+1)$ & $\vartheta>\omega(\alpha+1)$ \\ \hline
\rule[-1em]{0em}{2.5em} $\alpha\leq1$ & (\cAa) (\cCd) & (\cBa) (\cAb) (\cAc)(\cAd) \\
\rule[-1em]{0em}{2.5em} $\alpha>1$ & (\cCa) (\cDd) & (\cDa) (\cBb) (\cBc) (\cBd) \\ \hline
\end{tabular}
\end{center}

\vspace{0.5em}

Explicitly, we have
\begin{align*}
f(\alpha,\vartheta)\leq&\,\Xi\left(\mathfrak{c}\cdot\alpha^{\epsilon/2}+\mathfrak{c}\cdot\alpha^{n-1+\epsilon}\right)
   \qquad\text{for $\alpha\leq1$ and $\vartheta\leq\omega(\alpha+1)$,} \\
f(\alpha,\vartheta)\leq&\,\Xi\left(\frac{\mathfrak{c}\cdot\alpha^{\epsilon/2}}{[\vartheta+\omega(1-\alpha)]^{2n}}
   +\frac{\mathfrak{c}\cdot\alpha^{n-1}}{[\vartheta+\omega(2-\alpha)]^{3n-1}} +\frac{\mathfrak{c}\cdot\alpha^{n+\epsilon}}{[\vartheta+\omega(2\alpha-1)]^{n+1}}\right) \\
&\qquad\text{for $\alpha\leq1$ and $\vartheta>\omega(\alpha+1)$,} \\
f(\alpha,\vartheta)\leq&\frac{\mathfrak{c}\,\Xi}{\alpha^{n+\epsilon}} \qquad\text{for $\alpha>1$ and $\vartheta\leq\omega(\alpha+1)$,} \\
f(\alpha,\vartheta)\leq&\,\Xi\left(\frac{\mathfrak{c}}{\alpha^{\epsilon/2}\,[\vartheta+\omega(\alpha-1)]^{n+\epsilon/2}}
   +\frac{\mathfrak{c}}{\alpha\,[\vartheta+\omega(2-\alpha)]^{2n-1}}+\frac{\mathfrak{c}\cdot\alpha^{n+\epsilon}}{[\vartheta+(2\alpha-1)]^{3n}}\right) \\
&\qquad\text{for $\alpha>1$ and $\vartheta>\omega(\alpha+1)$,}
\end{align*}
and hence, $f$ is an $\mathcal{L}^1$--integrable function over $\mathbb{K}$. Since the value of the integral depends linearly on the constant~$\Xi$, it can
be arbitrarily small. \hfill $\Box$\\

\begin{bfseries}Remark. \end{bfseries}In~\cite{IH10} the estimation~$ (\cDa)$ is not sufficient for the $\mathcal L^1$--convergence of function~$f$.

\section{Density results for wavelet frames}\label{sec:density}

Wavelet frames described in Theorem~\ref{thm:Pframes} are semi--regular, i.e., for each scale a discretization of the position is performed. As a next step, we prove that discrete wavelet frames exist if only the set of sampling points is dense enough with respect to scale and position simultaneously. For a precise formulation we introduce the notion of hyperbolic density of a grid.

\begin{df}\label{def:hyperbolic_density}
We say a grid of points inside the unit ball~$\mathbb{B}$ is \emph{of density~$\rho$} if any ball inside~$\mathbb{B}$ with radius~$\rho$ with respect to the metrics
\begin{align*}
d\zeta_h:=&\left\|\frac{2}{1-r^2}\,(dr,\,h\,d\psi)\right\|
\end{align*}
(with~$\psi$ denoting the geodesic distance of points on the sphere and $h$ -- a positive constant) contains at least one grid point. We say a grid $\Lambda\subseteq\mathbb R_+\times\mathcal S^n$ is of \emph{hyperbolic density}~$\rho$ if $\left\{(e^{-b},y):\,(b,y)\in\Lambda\right\}\subseteq\mathbb B$ is a grid of density~$\rho$.
\end{df}

Now, we can formulate the result as follows.

\begin{thm}
\label{thm:Pframesn3dens} Let~$\{g_a:\,a\in\mathbb{R}_+\}$ be a wavelet family satisfying~\eqref{eq:pseudo_generating} with~$\gamma$ such that $\int_0^\infty\left||\gamma^2|^\prime(t)\right|\,dt<\infty$ and the kernel~$\Pi$ satisfying
$$
\left.
\begin{array}{c}
|\Pi(x,a;y,b)| \\[0.5em]
\left|(a+b)\,\nabla_{\!\ast}\Pi(x,a;y,b)\right| \\[0.5em]
\left|a\frac{\partial}{\partial a}\,\Pi(x,a;y,b)\right|
\end{array}
\right\} \leq(ab)^{n+\epsilon}\cdot
\begin{cases}
\frac{\mathfrak{c}}{(a+b)^{3n+2\epsilon}}, & \hspace{-0.7em}_{\angle(x,y)\leq
\omega[a+(2-\tilde\epsilon)b],} \\
\frac{\mathfrak{c}}{\angle(x,y)^{3n+2\epsilon}}, & \hspace{-0.7em}
_{\angle(x,y)>\omega(a+\tilde\epsilon b)},
\end{cases}
$$
for $a,b\leq b_0$ and for some positive constants~$\mathfrak{c}$, $\omega$, $\epsilon$, and $\tilde\epsilon<1/2$. Then there exists a constant~$\rho$ such that for any grid~$\Lambda$ of hyperbolic density~$\rho$ the family $\{g_{b,y}:\,(b,y)\in\Lambda\}$ is a weighted frame for~$\mathcal{L}^2(\mathcal S^n)$.
\end{thm}

\begin{bfseries} Proof. \end{bfseries}We show that the grid has the following property: There exist constants~$\delta$, and~$\Xi$ and a decreasing sequence of scales $\mathcal{B}=(b_j)_{j\in\mathbb{N}_0}$ such that $b_0\geq-\log\rho$ and the ratio $b_j/b_{j+1}$ is uniformly bounded from above by~$1+\rho\cdot\delta$. At each scale \mbox{$b=b_j$}, there is a measurable partition of~$\mathcal S^n$ $\mathcal{P}_b=\{\mathcal{O}_k^{(b)}:\,k=1,2,\dots,K_b\}$ into simply connected sets such that the diameter of each set (measured in geodesic distance) is not larger than~$\rho\cdot\Xi b$. In any of the sets $(b_{j+1},b_j]\times\mathcal{O}_k^{(b_j)}$ there is at least one point of the grid.

The proof of this statement is analogous to the proof of \cite[Lemma~8]{IH10}, but the partition of the sphere which we have to choose for each radius~$r_j$ has the property that each of the sets~$\mathcal O_{jk}$ has a diameter not larger than $2\sqrt{n(n+1)}$ and inradius (i.e., the diameter of the inscribed spherical ball) larger than~$\rho$. The existence of such a partition is proven in Lemma~\ref{lem:partition}.

On the other hand, it is shown in Section~\ref{sec:semicontinuous_frames} that for any $\epsilon>0$ and $\delta>0$ there exists a constant~$\rho$ such that for any sequence~$\mathcal{B}=(b_j)_{j\in\mathbb{N}_0}$ with $b_0\geq-\log\rho$ and $1<b_{j}/b_{j+1}<1+\rho\cdot\delta$ the family \mbox{$\{g_{b_j,x},\,b_j\in\mathcal{B},\,x\in\mathcal S^n\}$} is a semi-continuous frame for~$\mathcal{L}^2(\mathcal S^n)$, satisfying the frame condition~\eqref{eq:framedef} with the prescribed~$\epsilon$.

Now, Theorem~\ref{thm:Pframesn3dens} can be proven in the same way as~\cite[Theorem~9]{IH10}.\hfill$\Box$

\begin{lem}\label{lem:partition}For any $d\leq1$ there exists a measurable partition of~$\mathcal S^n$ into simply connected sets such that the diameter of each set (measured in geodesic distance) is not larger than~$d$ and the radius of the inscribed spherical ball is larger than $\frac{d}{2\sqrt{n(n+1)}}$.
\end{lem}

\begin{bfseries}Proof.\end{bfseries} Consider a $n+1$--dimensional cube inscribed in the $n$--sphere. Its side length is equal to $\frac{2}{\sqrt{n+1}}$ and the diameter of each of its $2(n+1)$ facets equals $2\sqrt\frac{n}{n+1}$. Suppose each of the facets is subdivided into $2^{nk}$ $n$--dimensional cubes and consider the central projection of these cubes onto the sphere. The diameter of each subset is not larger than
\begin{equation}\label{eq:diameter_partition}
\arctan\frac{\frac{2}{2^k}\sqrt\frac{n}{n+1}}{\frac{1}{\sqrt{n+1}}}=\arctan\frac{\sqrt n}{2^{k-1}}\leq\frac{\sqrt n}{2^{k-1}},
\end{equation}
The numerator in~\eqref{eq:diameter_partition} is the diameter of each sub--cube, and the denominator is the closest distance of a cube from the origin of the axes. Further, each of these subsets contains a ball with a diameter equal to the side length of the original sub--cube, i.e., $\frac{1}{2^{k-1}\sqrt{n+1}}$. Consequently, for each~$d$ small enough, we can find a~$k$ such that
$$
\frac{\sqrt n}{2^{k-1}}\leq d\leq\frac{\sqrt n}{2^{k-2}}
$$
and for the $k^{\text{th}}$--level partition the inradius of each subset is greater than or equal to
$$
\frac{1}{2^{k-1}\sqrt{n+1}}\geq\frac{d}{2\sqrt{n(n+1)}}.
$$
\hfill$\Box$

\begin{bfseries}Remark.\end{bfseries} In the proof of~\cite[Lemma~7]{IH10} another partition for the two--dimensional sphere is proposed and the quotient of the radii is estimated more generously.

\section{Poisson wavelet frames}\label{sec:Poisson_frames}

In this section we show that Poisson wavelets satisfy the conditions of Theorem~\ref{thm:Pframesn3dens}, i.e., they yield discrete frames. In order to do it, we need to prove that the kernel~$\Pi^m$ is localized according to~\eqref{eq:extkern3}. Theorem~\ref{thm:scaling} and Corollary~\ref{cor:proporcja} will be used to characterize the localization of the kernel. For its gradient, we utilize the fact that the wavelets are rotation invariant, and hence, the modulus of their surface gradient is equal to the modulus of their derivative with respect to the first spherical variable~$\theta_1$. In a similar manner as in~\cite{IH10}, we derive estimations for~$|\nabla_{\!\ast} g_a^m|$. The proofs are analogous and we omit them here.

\begin{lem}
\label{lem:rop} Let $Q_m$, $m\in\mathbb{N}$, be a sequence of polynomials in two variables satisfying the recursion
$$
Q_{m+1}( r,t)=A_m( r,t)\cdot Q_m( r,t)+B( r,t)\cdot\frac{\partial}{\partial  r}Q_m( r,t)
$$
with
$$
A_m( r,t)=1-(\alpha+1) r^2+\alpha  rt\quad\text{for some positive }\alpha
$$
and
$$
B( r,t)=(1+ r^2-2 rt)\, r,
$$
and such that
\begin{equation}  \label{eq:assQ1}
Q_1(1,1)=0,\quad\text{and}\quad\left.\frac{\partial}{\partial  r}\,Q_1( r,1)\right|_{ r=1}\ne0.
\end{equation}
Then the polynomial $Q_m(1,\cdot)$, $m\geq2$, has an $\left[(m+1)/2\right]$--fold root in~$1$.
\end{lem}

\begin{lem}
\label{lem:envf} Let~$\{f_m\}$ be a family of functions over $(0,1)\times[0,\pi]$ given by
\begin{align*}
f_1( r,\theta)&=\frac{ r\,\sin\theta\,Q_1( r,\cos\theta)}{(1+ r^2-2 r\cos\theta)^{2+\lambda}}, \\
f_{m+1}( r,\theta)&= r\,\frac{\partial}{\partial  r}\,f_m( r,\theta),
\end{align*}
where~$Q_1$ is a polynomial satisfying~\eqref{eq:assQ1} and $\lambda$ is a positive integer or half--integer. Then, for any $k\geq2[m/2]+2\lambda+1$ there exists a constant~$\mathfrak{c}$ such that
\begin{equation}  \label{eq:envf}
|f_m( r,\theta)|\leq\mathfrak{c}\cdot\frac{ r}{\theta^k},\quad\theta\in(0,\pi],
\end{equation}
uniformly in~$ r$. For $m\geq2$, the number $2[m/2]+2\lambda+1$ is the smallest possible exponent~$k$. If~$Q_1(1,y)$ has a simple root in~$1$, then $2\lambda+1$ is
the smallest possible exponent~$k$ on the right--hand--side of~\eqref{eq:envf} for $m=1$.
\end{lem}

\begin{bfseries}Remark.\end{bfseries} In the proof of \cite[Lemma 5]{IH10}, the factor $\theta^k$ is missing in the definition of~$F(0,\theta)$.

\begin{thm}
Let
\begin{equation*}
\Psi_{ r\hat e}^m=\frac{1}{\Sigma_n}\sum_{l=0}^\infty l^m r^l\,C_l^\lambda,\qquad m\in\mathbb{N}_0,
\end{equation*}
be a field on the sphere caused by the multipole (monopole for $m=0$) $\mu=( r\partial_ r)^m\delta_{ r\hat{e}}$. For any $k\geq2\left[\frac{m}{2}\right]+n$ there exists a constant~$\mathfrak{c}$ such that
$$
\left|\frac{d}{d\theta}\,\Psi_{ r\hat e}^m(\cos\theta)\right|\leq\mathfrak{c}\cdot\frac{ r}{\theta^k},\quad\theta\in(0,\pi],
$$
uniformly in~$ r$. $2\left[\frac{m}{2}\right]+n$ is the smallest possible exponent on the right--hand--side of this inequality.
\end{thm}

\begin{cor}\label{cor:localization} For any $k\geq2\left[\frac{m+1}{2}\right]+n$ there exists a constant~$\mathfrak{c}$ such that
$$
\left|\frac{d}{d\theta}\,g_a^m(\cos\theta)\right|\leq\mathfrak{c}\cdot\frac{a^m\,e^{-a}}{\theta^k},\quad\theta\in(0,\pi],
$$
uniformly in~$ r$. $2\left[\frac{m+1}{2}\right]+n$ is the smallest possible exponent on the right--hand--side of this inequality.
\end{cor}

\begin{thm} Let~$g_a^m$ be a Poisson wavelet of order~$m$. Then there exists a constant~$\mathfrak{c}$ such that
$$
\left|a^{n+1}\frac{d}{d\theta}\,g_a^m\left(\cos(\theta)\right)\right|_{\theta=a\theta}\leq\frac{\mathfrak{c}\cdot e^{-a}}{\theta^{m+n+1}},\quad\theta\in\left(0,\frac{\pi}{a}\right],
$$
uniformly in~$a$, and $m+n+1$ is the largest possible exponent.
\end{thm}

\begin{bfseries}Remark.\end{bfseries} In~\cite[Theorem~8]{IH10} the argument of~$g_a^d$ should be equal to~$\theta$, and $a\theta$ is substituted after derivation.

\begin{cor}\label{cor:proporcja} The functions $(a,\theta)\mapsto a^{n+1}\frac{d}{d\theta}\,g_a^m(\cos\theta)$ are bounded by $\mathfrak c\cdot e^{-a}$ uniformly in~$\theta$. $n+1$ is the smallest possible exponent in this inequality.
\end{cor}

Having these estimations, we may now come to the main theorem of this section.

\begin{thm}
\label{thm:framesn3} Let~$\{g_a^m\}$, $a\in\mathbb{R}_+$, be a Poisson wavelet family of order $m\geq n+1$. Then, there exist constants~$\delta$ and$\rho$ such
that for any grid~$\Lambda$ of a type $(\delta,\Xi)$ with $\Xi\leq\rho$ the family $\{g_{b,y}^m:\,(b,y)\in\Lambda\}$ is a weighted frame for $\mathcal{L}^2(\mathcal S^n)$.
\end{thm}

Although the proof is analogous to the proof of~\cite[Corollary~4]{IH10}, we repeat its first part in order to clarify some inaccuracies.

\noindent\emph{Proof. }The semi--frame condition for a Poisson wavelet family is verified in Section~\ref{sec:semicontinuous_frames}. It is to be checked that the estimations on the kernel and its gradient hold. The kernel~$\Pi^m$ is given by
\begin{equation*}
\Pi^m(a,x;b,y)=\mathfrak c\cdot\frac{(ab)^m}{(a+b)^{2m}}\,g_{a+b}^{2m}\left(\angle(x,y)\right),
\end{equation*}
compare formula~\eqref{eq:kernel}, and for the wavelet we have the estimation
\begin{equation*}
\left|g_{a+b}^{2m}\left(\angle(x,y)\right)\right|\leq\mathfrak{c}\cdot\frac{(a+b)^{2m}}{\angle(x,y)^{2m+n}}
\end{equation*}
uniformly in $\angle(x,y)$, $a$, $b$, and
\begin{equation*}
\left|g_{a+b}^{2m}\left(\angle(x,y)\right)\right|\leq\frac{\mathfrak{c}}{(a+b)^n}
\end{equation*}
cf.~Theorem~\ref{thm:scaling} and Corollary~\ref{cor:proporcja}. Therefore, we have
$$
\left|\Pi^m(a,x;b,y)\right|=\mathfrak c\cdot\frac{(ab)^m}{\angle(x,y)^{2m+n}}
   \qquad\text{and}\qquad\left|\Pi^m(a,x;b,y)\right|=\mathfrak c\cdot\frac{(ab)^m}{(a+b)^{2m+n}}
$$
Further, since
\begin{equation*}
\frac{(ab)^{m-n-\epsilon}}{\theta^{2(m-n-\epsilon)}}\leq\mathfrak{c}
\end{equation*}
for $\theta\geq\lambda(a+\epsilon b)$ and $\epsilon<1$ and
\begin{equation*}
\frac{(ab)^{m-n-\epsilon}}{(a+b)^{2(m-n-\epsilon)}}\leq1
\end{equation*}
for $\epsilon<1$, the inequalities~\eqref{eq:extkern3} are satisfied for the kernel. Estimations for the gradient are obtained in an analogous way. \hfill $\Box$\\

\begin{bfseries}Remark. \end{bfseries}Consider a grid of points inside~$\mathbb B$ corresponding to the grid~$\Lambda$, compare Definition~\ref{def:hyperbolic_density}. For Poisson wavelets it is a grid of their sources.\\

In order to find an appropriate weight function, we need to choose a set of scales~$\mathcal B$ and for each scale a partition of~$\mathcal S^n$ $\mathcal P_b$ such that in each of the sets $(b_{j+1},b_j]\times\mathcal O_k^{(b_j)}$ there is at least one point $y\in\Lambda$. The weight corresponding to such a point is the measure of the set $(b_{j+1},b_j]\times\mathcal O_k^{(b_j)}$ divided by the number of points in it, compare the proof of \cite[Theorem~9]{IH10}.\\

We can also consider a generalization of Poisson multipole wavelets to Poisson wavelets of a non--integer order,
$$
g_a^\nu=\sum_{l=0}^\infty(al)^\nu\,e^{-al}\,\mathcal K_l^\lambda,\quad\nu\in\mathbb R_+.
$$
It is straightforward to verify that for a half--integer~$\nu$ these functions satisfy the assumptions of Theorem~\ref{thm:Pframesn3dens} if $\nu>n$. Hence, if the hyperbolic density of~$\Lambda$ is large enough, the set $\{g_{b,y}^\nu,\,(b,y)\in\Lambda\}$ constitutes a weighted frame for $\mathcal L^2(\mathcal S^n)$.\\

\section{Frames of Mexican needlets}\label{sec:Mexican_needlets}

Recently a new wavelet construction for compact manifolds~\cite{GM09a} and its discretization to frames~\cite{GM09b} have been presented. The wavelets are kernels of the convolution operator $f(t\Delta^{\!\ast})$, where $0\ne f\in\mathcal S(\mathbb R_+)$, $f(0)\ne0$, and~$\Delta^{\!\ast}$ denotes the Laplace--Beltrami operator on a manifold. In the case of the sphere this leads to zonal wavelets of the form
\begin{equation}\label{eq:Mneedlets}
K_a(\hat e,y)=\frac{1}{\Sigma_n}\sum_{l=0}^\infty f\left(a^2l(l+2\lambda)\right)\frac{\lambda+l}{\lambda}\,C_l^\lambda(y).
\end{equation}
For $f(s)=f_r(s)=s^re^{-s}$, $r\in\mathbb N$ they are called Mexican needlets. We have proved in~\cite[Section~5]{IIN14CWT} that they satisfy (up to a constant) the conditions of Definition~\ref{def:sphwavelet} for $\alpha(a)=\frac{1}{a}$.

Consider recalled kernels~$K_{\sqrt a}$. Note that the measures $\frac{d\sqrt a}{\sqrt a}$ and $\frac{da}{a}$ differ only by a multiplicative constant. The kernel of the wavelet transform is (up to a constant) equal to
\begin{equation*}\begin{split}
\widetilde\Pi(a,x;b,y)&=\frac{1}{\Sigma_n}\sum_{l=0}^\infty(ab)^r\left[l(l+2\lambda)\right]^{2r}e^{-(a+b)l(l+2\lambda)}\mathcal K_l^\lambda(x\cdot y)\\
&=\frac{(ab)^r}{(a+b)^{2r}}\,\widetilde K_{a+b}^{2r}(x,y).
\end{split}\end{equation*}
Formula~(8) from~\cite{GM09a} yields the following estimation
\begin{equation}\label{eq:est_Mneedlets}
\left|\frac{(ab)^r}{(a+b)^{2r}}\,\widetilde K_{a+b}^{2r}(x,y)\right|\leq\frac{\mathfrak c\,(ab)^r}{(a+b)^{2r+\frac{n-N}{2}}\angle(x,y)^N}.
\end{equation}
Suppose, one wants to obtain the estimation~\eqref{eq:extkern3} of the modulus of the kernel for big angles~$\angle(x,y)$. Since $\angle(x,y)$ is not bounded from below by a constant, $N$ has to be chosen to be less than or equal to $3n+2\epsilon$. Further, since $\angle(x,y)$ is not bounded from above by $\mathfrak c(a+b)$ for any~$\mathfrak c$, the exponent~$N$ cannot be less than~$3n+2\epsilon$. For this choice of~$N$ one obtains from~\eqref{eq:est_Mneedlets}
\begin{equation}\label{eq:est_Mneedlets_kernel}
\left|\widetilde\Pi(a,x;b,y)\right|\leq\frac{\mathfrak d\,(ab)^{n+\epsilon}}{\angle(x,y)^{3n+2\epsilon}}\cdot\frac{(ab)^{r-n-\epsilon}}{(a+b)^{2r-2n-2\epsilon}(a+b)^{n+\epsilon}}
\end{equation}
for $\angle(x,y)>\omega(a+\tilde\epsilon b)$. The second factor on the right--hand--side of~\eqref{eq:est_Mneedlets_kernel} cannot be estimated from above. Consequently, Theorem~\ref{thm:Pframes} does not apply to Mexican needlets.

\begin{bfseries}Remark.\end{bfseries} The discretization method could be adapted to Mexican needlets, but the estimations appearing in the proof of Theorem~\ref{thm:Pframes} would need to be repeated. Another possibility is to investigate whether Mexican needlets satisfy the estimations~\eqref{eq:extkern3}; inequality~\eqref{eq:est_Mneedlets_kernel} is a consequence of \cite[Lemma~4.1]{GM09a} which is proved in a much more general setting and with quite sophisticated methods. On the other hand, the discretization method presented here could also be applied to other homogeneous manifolds if appropriate assumptions on the kernel are made. An example of its application to Gabor transform over~$\mathbb R$ can be found in~\cite{IH10}.

A frame construction for Mexican needlets is described in~\cite{GM09b}. Here we want to discuss the similarities and differences in two approaches.

In both cases, first, scale discretization and then position discretization of the image of the wavelet transform is performed. For the latter one writes~$\mathcal S^n$ as a disjoint union of measurable sets and picks a point from each of the sets in which a frame vector is evaluated. Further, in both cases one deals (in principle) with weighted frames; in~\cite{GM09b} the weight (equal to the square root of the measure of the set) is included as a factor in each of the frame vectors. The essence of the proofs is an estimation of the error that occurs by replacing an integral over the sphere by its discretization.

There are several advantages of the method presented in~\cite{GM09b}. First, it is valid for any smooth compact oriented Riemannian manifold without a boundary. Second, a convergence speed to nearly tight frames is computed. Next, scales summation is truncated so that one obtains a finite sum in the frame inequality.

On the other hand, the auxiliary function~$G_{j,k}$ in the proof of \cite[Theorem~2.4]{GM09b} contains a set measure in the numerator, hence, the sets in the partition of~$\mathcal S^n$ cannot be too small. Further, the present method relying on a comparison of discrete and continuous convolution of kernels is applicable also when the gradient of a single kernel is not bounded (although Theorem~\ref{thm:Pframes} would require a slight reformulation). An example of such a kernel (for the Gabor transform) can be found in~\cite[Section~4]{IH13}. In our approach, scale discretization is more flexible, and scale perturbation (for single sampling points) after position discretization is studied.

\end{document}